\documentclass[12pt,reqno]{amsart}
\usepackage{verbatim}
\usepackage{bm}
\usepackage{amsfonts}
\usepackage{amsmath}
\usepackage{amssymb}
\newtheorem{theorem}{Theorem}

\newtheorem{remark}{Remark}

\newfont{\bb}{msbm10 at 12pt}

\def\qed{\hfill{Q.E.D.}\smallskip}

\newcommand{\ls}{\setlength{\baselineskip}{12pt}
                 \setlength{\parskip}{3mm}}

\newcommand{\mysection}[1]{\section{#1}\setcounter{equation}{0}}

\newcommand{\bal}{\begin{align}}      \newcommand{\eal}{\end{align}}
\newcommand{\ba}{\begin{array}}      \newcommand{\ea}{\end{array}}
\newcommand{\bc}{\begin{center}}     \newcommand{\ec}{\end{center}}
\newcommand{\be}{\begin{enumerate}}  \newcommand{\ee}{\end{enumerate}}
\newcommand{\beq}{\begin{eqnarray}}  \newcommand{\eeq}{\end{eqnarray}}
\newcommand{\beQ}{\begin{eqnarray*}} \newcommand{\eeQ}{\end{eqnarray*}}
\newcommand{\bi}{\begin{itemize}}    \newcommand{\ei}{\end{itemize}}
\newcommand{\bt}{\begin{tabular}}    \newcommand{\et}{\end{tabular}}
\newcommand{\bdm}{\begin{displaymath}} \newcommand{\edm}{\end{displaymath}}

\begin{document}

\title{ Lower bounds for the eigenvalue estimates of the  submanifold Dirac operator }

{\address{Yongfa Chen
 School of Mathematics and Statistics,
Central China Normal University, Wuhan 430079, P.R.China.}

\author{Yongfa Chen}}

\thanks{ }

\email{yfchen@mail.ccnu.edu.cn}

\begin{abstract}
We get optimal lower bounds for the eigenvalues of  the submanifold Dirac operator on locally reducible Riemannian manifolds
 in terms of intrinsic and extrinsic expressions. The limiting-cases are also studied.
As a corollary, one gets several known results in this direction.
\end{abstract}

\keywords{  Dirac operator,  eigenvalue, mean curvature,
scalar curvature}

\subjclass[2008]{}

\maketitle

\mysection{\textbf{Introduction}} \ls
It is well known  that  the spectrum of
the Dirac operator on closed  spin manifolds detects
subtle information on the geometry and the topology of such
manifolds (see \cite{LM}).
A fundamental tool to get estimates for eigenvalues of  the basic Dirac operator $D$ acting on spinors
is the Schr\"{o}dinger-Lichnerowicz  formula
\beq \label{S-L1}
D^2=\nabla^*\nabla+\frac{1}{4}{R}^M,
\eeq
where $\nabla^*$ is the formal adjoint of $\nabla$  with respect to the natural Hermitian
inner\\
 product
 on spinor bundle $\Sigma {M^m}$ and  $R^M$ stands for  the
scalar curvature of  the closed spin manifold $(M^m,g)$.

   The problem of finding optimal lower bounds for the  eigenvalues of the  Dirac operator on closed manifolds was  for the first time considered in 1980 by Friedrich. Using the Schr\"{o}dinger-Lichnerowicz formula and a modified spin connection, he proved the following sharp inequality:
 \beq\label{Friedrich}
 \lambda^2\geqslant c_m \min_M R^M,
\eeq
where $c_m=\frac{m}{4(m-1)}. $
The case of equality in (\ref{Friedrich}) occurs iff  $(M^m,g)$  admits a nontrivial spinor field
  $\psi$ called   a real \emph{Killing spinor},
satisfying
 the following overdetermined elliptic equation
\beq
\nabla_X\psi=-\frac{\lambda}{m}X\cdot \psi,
\eeq
where  $\lambda\in \mathbb{R}, \forall X\in \Gamma(TM)$ and the dot ``." indicates the Clifford multiplication.
The manifold  must be a locally irreducible Einstein manifold. Note complete simply-connected Riemannian spin manifolds $(M^m, g)$ carrying a
non-zero space of real Killing spinors have been completely classified by B\"ar\cite{Bar93}.
The limiting manifold must be either a stand $n$-sphere, an Einstein-Sasaki manifold, a $6$-dimensional nearly K\"ahler manifold or a $7$-dimensional manifold with $3$-form $\nabla\eta=\ast\eta.$

 The dimension dependent coefficient $c_m=\frac{m}{4(m-1)} $ in the estimate can be improved if one  imposes geometric assumptions on the  metric.
Kirchberg \cite{Kirchberg86,Kirchberg90} showed that  for K\"ahler metrics  $c_m$ can be replaced by   $\frac{m+2}{4m} $ if the complex dimension $\frac{m}{2} $  is odd, and  by
$\frac{m}{4(m-2)} $ if  $\frac{m}{2} $  is even.
Alexandrov, Grantcharov, and Ivanov \cite{A} showed that  if there  exists  a non-trivial parallel 1-form on $M^m$, then $c_m$ can be replaced by $c_{m-1}$. Later, Moroianu and Ornea \cite {Moroianu} weakened the assumption on the $1$-form from parallel to harmonic with constant length. Note the condition that the norm of the 1-form being constant is essential, in the sense that the topological constraint alone (the existence of a non-trivial harmonic $1$-form) does not allow any improvement of Friedrich's inequality (see \cite{Bar04}).
The generalization of \cite{A} to locally reducible Riemannian manifolds
was achieved by Alexandrov \cite{A1}, extending earlier work by  Kim
\cite{K}.

 On the other hand, motivated by the using of the hypersurface Dirac operator in the proof of
 the positive energy  conjecture by Witten, Zhang first investigated the corresponding eigenvalue
 estimates problem. Optimal lower bounds for the  hypersurface Dirac operator in terms of the scalar curvature, the mean curvature and energy-momentum tensor were obtained. Follow-up related results can be seen in \cite{Zhang, Hijazi1,Hijazi2, Morel, Ginoux}.

In this paper, based on  the Schr\"{o}dinger-Lichnerowicz  formula, by defining  appropriate modified connections, we shall  get  optimal lower bounds for the eigenvalues of  the submanifold Dirac operator on any locally reducible Riemannian manifolds in terms of the scalar curvature
 as well as a normal curvature term which only appears in codimension greater than one.
 These results can be also translated  in an intrinsic way for a twisted Dirac-Schr\"odinger operator. As a corollary, one recovers  several classical results in \cite{A1,A,Ginoux,Hijazi2,Moroianu}. The discussion of the limiting cases of these estimates give rise to two new field equations generalizing the
Killing equation.

The remainder of the article is organized as follows. Section 2 and 3 describe some geometric conventions and preliminaries about ``partial" Dirac operator, the $\beta$-twist $D_\beta$ and
the submanifold Dirac operator on locally reducible Riemannian manifolds. The main result and its proof are given in section 4.
In section 5, we consider the intrinsic estimates  for the twisted Dirac-Schr\"odinger operator.
In the final section, we  shall extend the our techniques to other generalized hypersurface Dirac operators which  appear in General Relativity.

\mysection{\textbf{On locally reducible Riemannian  manifolds}}

Let $M^m$ be a closed Riemannian spin manifold with positive scalar $R^M.$ Suppose
$TM=T_1\oplus\cdots\oplus T_k $ is orthogonal sum, where $T_i$ are parallel distributions of dimension $m_a,a=1,\cdots,k,$ and $m_1\geq m_2\geq\cdots\geq m_k.$ One important
consequence of the parallelism of $T_a$ is that $R(X,Y)=0$ whenever $X\in T_a, Y\in T_b$ with
$a\neq b.$
Then one can define
a locally decomposable
Riemannian structure $\beta$ as follows
\beq
\beta\mid_{T_1}=\textrm{Id},\ \  \beta\mid_{T_1^{\perp}}=-\textrm{Id}.
\eeq
Suppose
$$\{e_1,\cdots,e_{m_1},e_{m_1+1}, \cdots,e_{m_1+m_2 },\cdots,e_{m_1+m_2+\cdots+m_{k-1}+1 },\cdots, e_m\}$$
is an adapted local orthonormal frame, i.e., such that
 $\{e_{m_1+m_2+\cdots+m_{a-1}+1},\cdots,\\
 e_{m_1+m_2+\cdots+m_a}\}$ spans the subbundle $T_a.$
 Let $I_a=\{m_{a-1}+1,\cdots,m_a\},$
then
\beq
D_a&=&\sum_{i\in I_a} e_i\cdot\nabla_i
\eeq
is the ``partial" Dirac operator of subbundle $T_a,$ which  is formally self-adjoint operator.
Hence the Dirac operator $D$ and the $\beta$-twist $D_\beta$ can be expressed as the following
\beq
D&=&\sum_{k=1}^me_k\cdot\nabla_k=D_1+D_2+\cdots+D_k\\
D_\beta&:=&\sum_{k=1}^m\beta(e_k)\cdot\nabla_k=D_1-(D_2+\cdots+D_k)
\eeq
Note for $a\neq b,$ one has
\beq\label{DaDb}
D_aD_b+D_aD_b=0,
\eeq
that is,
\beq\label{Dsquare}
D_\beta^2=D^2=\sum_{a=1}^kD^2_a.
\eeq
In addition, we also have
\beQ
DD_\beta+D_\beta D=2(D_1^2-D_2^2-\cdots-D_k^2)
\eeQ

Let $R_a^M$ be  the ``scalar curvature" of $T_a,$ i.e.,
$$R_a^M=\sum_{s,t\in I_a}\langle R(e_s,e_t)e_t, e_s\rangle.$$
Hence the scalar curvature $R^M$ of $M^m$ is $R^M=\sum_{a=1}^kR_a^M.$

\mysection{\textbf{The submanifold Dirac operator}} \ls

Now let $(\overline{M}^{m+n},\bar{g})$ be an $(m+n)$-dimensional Riemannian spin manifold  and $M^m$ be an $m$-dimensional submanifold
in ${\overline{M}}^{m+n}$ with its induced  Riemannian structure. Assume that two manifolds are equipped with a spin structure
so that  there exists unique spin  structure on the normal bundle of $M^m$ such that the sum of the spin structures of the tangent bundle $TM^n$ and of the normal bundle $NM^n$ of $M^n$ is exactly the spin structure of $\overline{M}^{m+n}$ restricted to $M^m.$
Note that  in particular $M^m,\overline{M}^{m+n}$ are oriented.
Denote $\Sigma \overline{M}$  the spinor bundle of $\overline{M}^{m+n},$ then $\mathbb{S}:={\Sigma \overline{M}|}_{M}$ is globally defined along  $M^m.$

  Denote  the Levi-Civita  connections of $M^m$ and ${\overline{M}}^{m+n}$ by
  ${\nabla}$ and $\overline{\nabla}$ respectively and denote by the same symbol their corresponding  lift to the
  spinor  bundle $\mathbb{S}.$ Consider the Dirac operator $D$ of $M^m$ defined by $\nabla$
   on $\mathbb{S}$  and $\overline{D}$ the submanifold Dirac Operator defined by $\overline{\nabla}$ on $\mathbb{S}.$
   It is  known that  there exists  a positive definite Hermitian metric on   $\mathbb{S}$
   which satisfies, for any $X\in \Gamma(T\overline{M}),$ any spinor fields $\phi,\varphi\in\Gamma(\mathbb{S}),$ the relation
\beQ
\langle X \cdot\phi,X \cdot\varphi\rangle=|X|^2\langle\phi,\varphi\rangle,
\eeQ
where $`` \cdot " $denotes Clifford  multiplication on $\overline{M}^{m+n}.$ This metric is globally defined
along $M^m.$ The connection is compatible with the metric $\langle \cdot,\cdot\rangle.$ Fix a point $p\in M^m$ and
a local orthonormal basis $\{e_\alpha\}$ of $T_p\overline{M}^{m+n}$ with $ \{e_A\}$  normal to $M^m$ and $\{e_i\}$ tangent to $M^m$ such that
for $1\leq i,j\leq m$
\beQ
(\nabla_ie_j)_p=0.
\eeQ
All the computations will be made in such charts.
 Then, at point $p,$ for $1\leq i,j\leq m$ and $m+1\leq A,B\leq m+n,$
\beq
{\overline{\nabla}}_ie_A&=&-\sum_{j} h_{Aij}e_j+{\nabla}_i^\bot e_A,\label{Gaussformula}\\ {\overline{\nabla}}_ie_j&=&\sum_{A}h_{Aij}e_A,\\
{\nabla}_i^\bot e_A&:=&\sum_{B}a_{iAB}e_B,
\eeq
where
\beQ
h_{Aij}=h_{Aji}=-\langle{\overline{\nabla}}_ie_A,e_j\rangle,\quad
a_{iAB}=-a_{iBA}=\langle{\overline{\nabla}}_ie_A,e_B\rangle
\eeQ
 are the components of the second fundamental form $h$ and normal connection ${\nabla}^\bot$ at $p$ respectively. The spinorial Gauss formula says that, for $1\leq i,j\leq m$ and $m+1\leq A\leq m+n,$
\beq\label{spinorialGaussformula}
{\overline{\nabla}}_i={\nabla}_i+\frac{1}{2}\sum_{j,A}h_{Aij} e_j\cdot e_A\cdot,
\eeq
which implies that the connection $\nabla$ is compatible with the metric $\langle \cdot,\cdot\rangle$ too, and also  for any spinor fields $\phi\in\Gamma(\mathbb{S}),$
\beq
\nabla_i(e_A\cdot\phi)
&\stackrel{(\ref{spinorialGaussformula})}{=}&{\overline{\nabla}}_i( e_A\cdot\phi)-\frac{1}{2}\sum_{j,B}h_{Bij} e_j\cdot e_B\cdot e_A\cdot\phi\nonumber\\
&=&{\overline{\nabla}}_i e_A\cdot\phi+ e_A\cdot{\overline{\nabla}}_i \phi+\frac{1}{2}\sum_{j,B}h_{Bij} e_j\cdot( e_A\cdot e_B+2\delta_{AB})\cdot\phi\nonumber\\
&\stackrel{(\ref{Gaussformula})}{=}&\left(-\sum_{j} h_{Aij}e_j+{\nabla}_i^\bot e_A\right)\cdot\phi+ e_A\cdot{\overline{\nabla}}_i \phi-e_A\cdot\frac{1}{2}\sum_{j,B}h_{Bij} e_j\cdot e_B\cdot\phi\nonumber\\
&{}&+\sum_{j}h_{Aij} e_j\cdot\phi\nonumber\\
&\stackrel{(\ref{spinorialGaussformula})}{=}&\nabla_i^\bot e_A\cdot\phi+e_A\cdot\nabla_i\phi.
\eeq

Note that  in the above orthonormal frame $\{e_i\}$ of $M^m,$ the  Dirac operator and the submanifold Dirac operator is defined as follow,
\beq
D=\sum_{i=1}^m e_i\cdot \nabla_i,\ \ \ \ \ \ \
{\overline{D}}=\sum_{i=1}^m e_i\cdot {\overline{\nabla}}_i,
\eeq
 respectively. Contrast to the basic Dirac operator, ${\overline{D}} $ is in general not self-adjoint w.r.t. positive definite Hermitian metric on $\mathbb{S}.$  In fact, it is easy to see that
 \beq
{\overline{D}}&=&D-\frac{1}{2}\sum_{A}H_Ae_A\cdot,\\
\overline{D}^{*}&=&D+\frac{1}{2}\sum_{A}H_Ae_A\cdot,
 \eeq
where $H_A:=\sum_{i}h_{Aii}$ is the component of the mean curvature vector field of $M^m.$  From now on, we always denote $\sum_AH_Ae_A$ by $\vec{H}.$

Recall $\omega_n=i^{[\frac{n+1}{2}]}e_{m+1}\cdots e_{m+n}$ is the complex volume element of the normal bundle $NM^n.$
There is an operator $``\omega_\bot"$ on $\mathbb{S}$ defined by
\beq
\omega_\bot:=\left\{
    \begin{array}{cccccc}
\omega_n\quad\quad  \textup{if}\   n\ \textup{is even}\\
-i\omega_n\quad \textup{if}\  n\  \textup{is odd}\\
\end{array}
    \right.
\eeq
(see \cite{Ginoux,Hijazi2}).
Then one can check the following
\beQ
\omega_\bot&=&(-1)^n,\\
\langle\omega_\bot\cdot\phi,\varphi \rangle&=&(-1)^n\langle\phi,\omega_\bot\cdot\varphi \rangle,\\
\omega_\bot\cdot e_A\cdot&=&(-1)^{n-1}\omega_\bot\cdot e_A\cdot
\eeQ
and
\beQ
\nabla_i(\omega_\bot\cdot\phi)&=&\omega_\bot\cdot\nabla_i\phi,\\
D(\omega_\bot\cdot\phi)&=&(-1)^n\omega_\bot\cdot D\phi.
\eeQ
Hence
\beq
\overline{D}^{*}=\omega_\bot\cdot \overline{D}\omega_\bot\cdot.
\eeq
One can also prove that
\beq
D_H&:=&(-1)^n\omega_\perp\cdot {\overline{D}}\\
&=&(-1)^n\omega_\bot\cdot D+\frac{1}{2}\vec{H}\cdot\omega_\bot\cdot
\eeq
is formally self-adjoint with respect to
the metric $\langle \cdot,\cdot\rangle$ and, hence, $D_H$ has real eigenvalues.  This  first-order Dirac operator
arises in apparent  horizons in the spinor  proof of the positive mass theorem for black holes.
Moreover,  we have the following Schr\"{o}dinger-Lichnerowicz  type formula
\beq\label{SSL}
{\overline{D}}^{*}{\overline{D}}\phi&=&D_H^2\phi\nonumber\\
&=&\left((-1)^n\omega_\bot\cdot D+\frac{1}{2}\vec{H}\cdot\omega_\bot\cdot\right)
\left((-1)^n\omega_\bot\cdot D\phi+\frac{1}{2}\vec{H}\cdot\omega_\bot\cdot\phi\right)\nonumber\\
&=&D^2\phi+\frac{(-1)^n}{2}\omega_\bot\cdot D(\vec{H}\cdot \omega_\bot\cdot\phi)
+\frac{1}{2}\vec{H}\cdot D\phi+\frac{1}{4}|\vec{H}|^2\phi\nonumber\\
&=&D^2\phi+\frac{1}{4}|\vec{H}|^2\phi-\frac{1}{2}D^\bot\vec{H}\cdot\phi+\vec{H}\cdot D\phi\nonumber\\
&=&\nabla^{*}\nabla\phi+\frac{1}{4}(R^M+{\mathfrak{R}}_{\phi}^N+|\vec{H}|^2)\phi
-\frac{1}{2}D^\bot\vec{H}\cdot\phi+\vec{H}\cdot D\phi,\quad
\eeq
here $D^\bot:=\sum_{i=1}^m e_i\cdot \nabla_i^\bot$
and for any spinor field $\phi\in\Gamma(\mathbb{S}),$  the real function
\beq
{\mathfrak{R}}_{\phi}^N:=-\frac{1}{2}\sum_{i,j,A,B}R_{ijAB}\langle e_i\cdot e_j\cdot e_A\cdot e_B \cdot\phi,\phi/|\phi|^2\rangle
\eeq
is defined on subset $M_{\phi}:=\{x\in M^m|\phi(x)\neq 0\},$ where ${R}_{ijAB}$ is the curvature tensor of the
normal bundle $NM^m.$

\mysection{\textbf{Estimates  for the eigenvalues of the submanifold Dirac operator}} \ls
In this section, we introduce some modified connections to get the lower bounds for
the eigenvalues of the operator $D_H.$

\begin{theorem}\label{theorem1}
Let $M^m\subset \overline{M}^{m+n}$ be a closed spin submanifold of dimension $m\geq 2$ whose normal bundle is also spin. Suppose $TM^m=T_1\oplus\cdots\oplus T_k,$ where $T_a$ are parallel and pairwise orthogonal distributions of dimension $m_a, a=1,\cdots,
k, $ and
$m_1\geq m_2\geq\cdots\geq m_k\geq0.$ Consider a non-trivial spinor field $\psi\in \Gamma(\mathbb{S})$ such that
$D_H\psi=\lambda_H\psi,$ then
\beq\label{main result H1}
\lambda_H^2\geq\frac{1}{4}\sup_q\inf_{M_\psi}\left(\frac{R^M+\mathfrak{R}_{\psi}^N}{m_1q^2-2q+1}
-\frac{(m_1-1) |\vec{H}|^2}{(1-m_1q)^2}\right),
\eeq
where $q$ is  some real function,  $q\neq \frac{1}{m_1}$ if $\vec{H}\neq 0.$
If $\lambda_H^2$ achieves its minimum and the normal bundle is flat,  the scalar curvature $R^M$ and $|\vec{H}|$ are both constants.
\end{theorem}

\proof
First, one can define
a locally decomposable
Riemannian structure $\beta$ as follows
\beQ
\beta\mid_{T_1}=\textrm{Id},\ \  \beta\mid_{T_1^{\perp}}=-\textrm{Id}.
\eeQ
 and  we also define the following modified connection
\beq\label{modicon}
T_i\phi
&=&\nabla_i\phi
+\frac{1}{2}(\beta(e_i)+e_i)\cdot\left(\frac{1}{2}p\vec{H}+q\lambda_H\omega_\perp\right)\cdot\phi+
\frac{1}{2}\nabla_{(\beta-\textrm{Id})(e_i)}\phi
\eeq
where $p, q$ are smooth real-valued functions that are specified later.
Then, a direct  computation gives
\beq
|T\psi|^2&=&|\nabla \psi|^2+(m_1pq-p-q)\mathfrak{R}e\langle \vec{H}\cdot \psi,D\psi\rangle+(m_1q^2-2q)|D_H\psi|^2\nonumber\\
&{}&+\frac{|\vec{H}|^2}{4}(m_1p^2+2q-2m_1pq)|\psi|^2\nonumber\\
&{}&+\mathfrak{R}e\left\langle\left(\frac{1}{2}p\vec{H}+q\lambda_H\omega_\perp\right)\cdot\psi,
(D-D_\beta)\psi\right\rangle -\sum_{s=m_1+1}^m|\nabla_s\psi|^2,
\eeq
where $D_H\psi=\lambda_H\psi,$ for a non-trivial spinor field $\psi.$

Since $\mathfrak{R}e\langle D^\perp\vec{H}\cdot\phi, \phi\rangle=0,$  integrating (\ref{SSL}) over $M^m$ yields, for any spinor $\phi$
\beq
\int_M|D_H\phi|^2
=\int_M|\nabla\phi|^2+\frac{1}{4}(R^M+\mathfrak{R}^N_{\phi}+|\vec{H}|^2)|\phi|^2
-\mathfrak{R}e\langle \vec{H}\cdot\phi,D\phi\rangle.
\eeq
Hence we obtain, for the eigenspinor $\psi$ of $D_H$
\beq\label{integral term}
&{}&\int_M|T\psi|^2+\sum_{s=m_1+1}^m|\nabla_s\psi|^2
-\mathfrak{R}e\left\langle\left(p\vec{H}/2+q\lambda_H\omega_\perp\right)\cdot\psi,
(D-D_\beta)\psi\right\rangle\nonumber\\
&=&\int_M\left[m_1pq-p-q+1\right]\mathfrak{R}e\langle \vec{H}\cdot \psi,D\psi\rangle+[m_1q^2-2q+1]\lambda_H^2|\psi|^2\nonumber\\
&{}&-\frac{1}{4}(R^M+\mathfrak{R}^N_{\psi})|\psi|^2
+\frac{|\vec{H}|^2}{4}\left[m_1p^2+2q-2m_1pq-1\right]|\psi|^2.
\eeq
If $\vec{H}\neq0, $ let $m_1pq-p-q+1=0,$ that is, $$p=\frac{1-q}{1-m_1q}.$$
Then
\beq\label{IRHS1}
&{}&\textrm{R.H.S. of (\ref{integral term})}\nonumber \\
&=&\int_M\left(m_1q^2-2q+1\right)
\left[\lambda_H^2-\frac{1}{4}\left(\frac{R^M+\mathfrak{R}^N_{\psi}}{m_1q^2-2q+1}
-\frac{(m_1-1) |\vec{H}|^2}{(1-m_1q)^2}\right)\right]|\psi|^2.\ \ \ \ \  \ \ \ \
\eeq

Now we turn to deal with the L.H.S of (\ref{integral term}).
First, observe that for any spinor field $\phi\in \Gamma(\mathbb{S})$
\beq\label{T}
\sum_{i=1}^me_i\cdot T_i\phi=\frac{1}{2}(D_\beta+D)\phi-m_1\left(\frac{1}{2}p\vec{H}+q\lambda_H\omega_\perp\right)\cdot\phi
\eeq
Hence by $D_\beta^2=D^2,$ we can deduce the relation, for any spinor $\phi$:
\beq\label{key}
&{}&\mathfrak{R}e\int_M\langle \sum_{i=1}^me_i\cdot T_i\phi, -\frac{1}{2}(D_\beta\phi-D\phi)\rangle\nonumber\\
&=&-\frac{m_1}{2}\int_M\mathfrak{R}e\left\langle \left(p\vec{H}/2+q\lambda_H\omega_\perp\right)\cdot\phi,D\phi-D_\beta\phi
\right\rangle.
\eeq
And also note for any spinor $\phi$ and $s>m_1, T_s\phi=0.$  Hence
\beQ
\mathfrak{R}e\left\langle \sum_{i=1}^me_i\cdot T_i\phi, \frac{1}{2}(D_\beta\phi-D\phi)\right\rangle
&=&\left|\sum_{i=1}^{m_1}e_i\cdot T_i\phi\right|\left|\frac{1}{2}(D_\beta\phi-D\phi)\right|\\
&\leq&\sqrt{m_1}\left| T\phi\right|\left|\frac{1}{2}(D_\beta\phi-D\phi)\right|\\
&\leq&\frac{\sqrt{m_1}}{2}\left(\varepsilon^{-1}| T\phi|^2+\frac{\varepsilon}{4}|D_\beta\phi-D\phi|^2\right),
\eeQ
where $\varepsilon$ is  to be a fixed   positive constant.

 Hence, plugging Equation (\ref{key}) into (\ref{integral term}) and  using the Cauchy-Schwarz inequality and the equality (\ref{Dsquare}), we have that
\beq\label{ILHS}
&{}&\textup{L.H.S. of}\ (\ref{integral term})\nonumber\\
&=&\int_M\left(1-\frac{1}{\varepsilon\sqrt{m_1}}\right)|T\psi|^2+\sum_{s=m_1+1}^m|\nabla_s\psi|^2
-\frac{\varepsilon}{\sqrt{m_1}}(|D_2\psi|^2+\cdots+|D_k\psi|^2)\nonumber\\
&\geq&\int_M\left(1-\frac{1}{\varepsilon\sqrt{m_1}}\right)|T\psi|^2
+\left(\frac{1}{m_2}-\frac{\varepsilon}{\sqrt{m_1}}\right)|D_2\psi|^2\nonumber\\
&&+\cdots+\left(\frac{1}{m_k}
-\frac{\varepsilon}{\sqrt{m_1}}\right)|D_k\psi|^2. \ \ \ \ \ \ \ \
\eeq
Now we take $\varepsilon=\frac{1}{\sqrt{m_2}},$
hence, (\ref{ILHS}),  combined with   (\ref{IRHS1}),  yields that
\beq\label{final}
0&{\leq}&\int_M\left(1-\sqrt{\frac{m_2}{m_1}}\right)|T\psi|^2
+\left(\frac{1}{m_2}-\frac{1}{\sqrt{m_1m_2}}\right)|D_2\psi|^2\nonumber\\
&&+\cdots+\left(\frac{1}{m_k}
-\frac{1}{\sqrt{m_1 m_2}}\right)|D_k\psi|^2\nonumber\\
&\leq&\int_M(m_1q^2-2q+1)\left[\lambda_H^2-\frac{1}{4}\left(\frac{R^M+\mathfrak{R}^N_{\psi}}{m_1q^2-2q+1}
-\frac{(m_1-1) |\vec{H}|^2}{(1-m_1q)^2}\right)\right]|\psi|^2.\quad\quad
\eeq
Therefore, the first part of the theorem follows.

If $\lambda_H^2$ achieves its minimum, then
\beq
\sum_{i\in I_1}e_i\cdot T_i\psi=\frac{1}{2}(D_\beta\psi-D\psi),
\eeq
 which, together with (\ref{T}), implies that,
\beq\label{D1}
D\psi=m_1\left(\frac{1}{2}p\vec{H}+q\lambda_H\omega_\perp\right)\cdot\psi.
\eeq
But $p=\frac{1-q}{1-m_1q}$ and
\beQ
D\psi&=&\left(\omega_\bot \cdot D_H+\frac{1}{2}\vec{H}\cdot\right)\psi\\
&=&\lambda_H\omega_\bot\cdot\psi+\frac{1}{2}\vec{H}\cdot\psi.
\eeQ
This implies that
\beq\label{Hcdot}
(m_1-1)\vec{H}\cdot\psi=2(1-m_1q)^2\lambda_H\omega_\perp\cdot\psi
\eeq
and moreover from (\ref{D1}),
\beq\label{D2}
D\psi=m_1\left(\frac{1}{2}p\vec{H}+q\lambda_H\omega_\perp\right)\cdot\psi=\frac{m_1}{2}\mathcal{H}\cdot\psi,
\eeq
where $\mathcal{H}:=\frac{1+m_1q^2-2q}{(1-m_1q)^2}\vec{H}.$

\textbf{Case 1}. If $m_1>m_2,$
the eigenspinor corresponding to  the smallest eigenvalue of  $D_H^2$ satisfies the following generalized Killing type equations
$$\nabla_{m_1+1}\psi=\nabla_{m_1+2}\psi=\cdots=\nabla_{m}\psi=0,$$
 as well as $T\psi=0,$ which is equivalent to, for $i=1,\cdots,m_1,$
\beq
\nabla_i\psi
=-e_i\cdot\left(\frac{1}{2}p\vec{H}+q\lambda_H\omega_\perp\right)\cdot\psi
=-\frac{1}{2}e_i\cdot\mathcal{H}\cdot\psi.
\eeq
Obviously, $d|\psi|^2=0. $ And for any $i,j=1,\cdots,m_1,$ $ s,t=m_1+1,\cdots,m, $
 one has $\nabla_i\nabla_s\psi=\nabla_s\nabla_t\psi=0,$
\beq
\nabla_s\nabla_i\psi=-\frac{1}{2}e_i\cdot\nabla_s^\bot\mathcal{H}\cdot\psi
\eeq
as well as
\beq
\nabla_j\nabla_i\psi
=-\frac{1}{2}e_i\cdot\nabla_j^\bot\mathcal{H}\cdot\psi+\frac{1}{4}|\mathcal{H}|^2e_i\cdot e_j\cdot\psi.
\eeq
If the normal bundle is flat, then
\beQ
0=-\frac{1}{2}Ric(e_s)\cdot\psi
&=&\sum_{j=1}^{m_1}e_j\cdot \nabla_{s}\nabla_j\psi\\
&=&\sum_{j=1}^{m_1}e_j\cdot\left(-\frac{1}{2}e_j\cdot\nabla_s^\bot\mathcal{H}\cdot \psi\right)=\frac{m_1}{2}{\nabla_s^\bot}\mathcal{H}\cdot \psi,
\eeQ
which implies that
\beq\label{sH}
\nabla_s^\bot\mathcal{H}=0,\ \textup{for} \  s>m_1.
\eeq
Moreover
\beq\label{Ric1}
-\frac{1}{2}Ric(e_i)\cdot\psi
&=&\sum_{j=1}^{m_1}e_j\cdot \mathcal{R}_{e_i,e_j}\psi\nonumber\\
&=&\frac{1}{2}(-D^\bot\mathcal{H}\cdot e_i\cdot\psi+m_1\nabla_i^\bot\mathcal{H}\cdot\psi)-\frac{m_1-1}{2}|\mathcal{H}|^2 e_i\cdot\psi.\quad\quad
\eeq
Using the fact $\sum_{i\in I_1} e_i \cdot D^\bot\mathcal{H}\cdot e_i=-(m_1-2)D^\bot\mathcal{H},$ one gets
\beq
\frac{1}{2}R^M=(m_1-1)D^\bot\mathcal{H}\cdot \psi+\frac{m_1(m_1-1)}{2}|\mathcal{H}|^2\psi.
\eeq
Take the inner product of the above equality with $\psi$ and compare its real and imaginary parts to obtain
\beq
D^\bot\mathcal{H}\cdot\psi=0,\quad\quad  R^M=m_1(m_1-1)|\mathcal{H}|^2.
\eeq
As a consequence, $R^M$ is a constant, due to (\ref{final}).

Furthermore, note
\beQ
D^\bot\mathcal{H}\cdot e_i\cdot\psi
&=&-\sum_{j=1}^me_j\cdot e_i\cdot\nabla_j^\bot\mathcal{H}\cdot\psi\\
&=&\sum_{j=1}^m(e_i\cdot e_j+2\delta_{ij})\cdot\nabla_j^\bot\mathcal{H}\cdot\psi\\
&=& e_i \cdot D^\bot\mathcal{H}\cdot\psi+2\nabla_i^\bot\mathcal{H}\cdot\psi=2\nabla_i^\bot\mathcal{H}\cdot\psi,
\eeQ
which yields
\beq\label{Ric2}
\frac{1}{2}Ric(e_i)\cdot\psi
&=&-\frac{m_1-2}{2}\nabla_i^\bot\mathcal{H}\cdot\psi+\frac{m_1-1}{2}|\mathcal{H}|^2 e_i\cdot\psi.
\eeq
Hence from (\ref{sH}) and (\ref{Ric2}), one gets $\nabla^\bot \mathcal{H}=0, $ if $m_1\neq 2$ and $Ric(e_i)=(m_1-1)|\mathcal{H}|^2 e_i,$ for $i=1,\ldots,m_1$ as well as
\beQ
\lambda_H^2=\frac{(m_1-1)^2}{4(1+m_1q^2-2q)^2}|\mathcal{H}|^2.
\eeQ

\textbf{Case 2}. If $m_1=m_2=\cdots=m_l<m_{l+1}\leq\cdots \leq m_k,$  then
\beq
D_{l+1}\psi=D_{l+2}\psi=\cdots= D_k\psi=0
\eeq
and we also have, for $i,j\in I_1,$ the two spinor fields are proportional, i.e.,
\beq
e_i\cdot T_i\psi=e_j\cdot T_j\psi
\eeq
which, in turn implies that
\beq\label{i}
e_i\cdot \nabla_i\psi=e_j\cdot \nabla_j\psi.
\eeq
And also for any $a\in \{2,\cdots,l\},$ and any $s,t\in I_a,$
\beq\label{s}
e_s\cdot \nabla_s\psi=e_t\cdot \nabla_t\psi.
\eeq
and $\alpha\in (I_1\cup I_2\cup\cdots\cup I_l)^c,$
\beq\label{alpha}
\nabla_\alpha\psi=0.
\eeq
Therefore, for $i\in I_1,$
\beq
\nabla_i\psi &=&T_i\psi
-e_i\cdot\left(\frac{1}{2}p\vec{H}+q\lambda_H\omega_\perp\right)\cdot\psi\nonumber\\
&=&-\frac{1}{2m_1}e_i\cdot (D_\beta\psi-D\psi)- \frac{1}{m_1}e_i\cdot D\psi\nonumber\\
&\stackrel{(\ref{D2})}{=}&- \frac{1}{4}e_i\cdot\mathcal{H}\cdot\psi-\frac{1}{2m_1}e_i\cdot D_\beta\psi,
\eeq
and while, for any $b\in \{2,\cdots,l\},s\in I_b,$
\beq
\nabla_s\psi
&=&-\frac{1}{m_1}e_s\cdot D_b\psi\nonumber\\
&=& \frac{1}{2m_1}e_s\cdot (D_\beta\psi-D\psi)+\frac{1}{m_1}\sum_{a=2,\neq b}^l e_s\cdot D_a \psi\nonumber\\
&=&- \frac{1}{4}e_s\cdot\mathcal{H}\cdot\psi+\frac{1}{2m_1}e_s\cdot D_\beta\psi
+\frac{1}{m_1}\sum_{a=2,\neq b}^l e_s\cdot D_a \psi.
\eeq
Using the facts, for $s,t\in I_b,$
$$\sum_{s\in I_b}e_s\cdot(e_s\cdot e_t-e_t\cdot e_s)=(-2m_1+2)e_t,$$
and
\beQ
\sum_{s\in I_b}e_s\cdot(e_s\cdot \nabla_t-e_t\cdot\nabla_s)
&=&-m_1\nabla_t+\sum_{s\in I_b}(e_t\cdot e_s\cdot\nabla_s+2\delta_{st}\nabla_s)\\
&=&(-m_1+2)\nabla_t+e_t\cdot D_b,
\eeQ
it follows
\beQ\label{Ric3}
-\frac{1}{2}Ric(e_t)\cdot\psi
&=&\sum_{s\in I_b}e_s\cdot\mathcal{R}_{e_t,e_s}\psi
\nonumber\\
&=&\sum_{s\in I_b}e_s\cdot[\nabla_t,\nabla_s]\psi\nonumber\\
&=&\frac{1}{4}(m_1\nabla_t^\bot\mathcal{H}-D_b^\bot\mathcal{H}\cdot e_t)\cdot\psi-\frac{m_1-1}{8}|\mathcal{H}|^2e_t\cdot\psi\nonumber\\
&{}&-\frac{m_1-1}{4c}e_t\cdot\mathcal{H}\cdot D_\beta\psi+\frac{1}{2m_1}\left[(2-m_1)\nabla_t(D_\beta\psi)+e_t\cdot D_b(D_\beta\psi)\right]\nonumber\\
&{}&-\frac{m_1-1}{2m_1}e_t\cdot \mathcal{H}\cdot \sum_{a=2,\neq b}^l  D_a \psi
+\frac{1}{m_1}\sum_{a=2,\neq b}^l  \left[(2-m_1)\nabla_t(D_a \psi)+e_t\cdot D_b(D_a \psi)\right].
\eeQ
Performing its Clifford multiplication by $e_t$ and summing over $t\in I_b,b\in\{2,\cdots,l\}$ yields,
\beq\label{Scalb}
\frac{1}{2}R_b^M\psi
&=&\frac{1}{2}\sum_{t\in I_b} e_t\cdot Ric(e_t)\cdot\psi\nonumber\\
&=&\frac{m_1-1}{2}D_b^\bot\mathcal{H}\cdot\psi+\frac{m_1(m_1-1)}{8}|\mathcal{H}|^2\psi
+\frac{m_1-1}{4}\mathcal{H}\cdot D_\beta\psi-\frac{m_1-1}{m_1}D_b(D_\beta\psi)\nonumber\\
&{}&+\frac{m_1-1}{2}\mathcal{H}\cdot \sum_{a=2,\neq b}^l  D_a \psi
-\frac{2(m_1-1)}{m_1}\sum_{a=2,\neq b}^l D_b(D_a \psi).
\eeq
On the other hand,
\beq\label{Scal1}
\frac{1}{2}R_1^M\psi
&=&-\frac{1}{2}\sum_{i\in I_1} e_i\cdot Ric(e_i)\cdot\psi\nonumber\\
&=&\frac{m_1-1}{2}D_1^\bot(\mathcal{H})\cdot\psi+\frac{m_1(m_1-1)}{8}|\mathcal{H}|^2\psi\nonumber\\
&{}&-\frac{m_1-1}{4}\mathcal{H}\cdot D_\beta\psi+\frac{m_1-1}{m_1}D_1(D_\beta\psi).
\eeq
Combining Eqs. (\ref{Scal1}) and (\ref{Scalb}), by (\ref{DaDb}) one gets
\beq
\frac{1}{2}R^M \psi
&=&\frac{1}{2}\sum_{a=1}^l R_a^M \psi\nonumber\\
&=&\frac{m_1-1}{2}\sum_{a=1}^lD_a^\bot\mathcal{H}\cdot\psi+l\cdot\frac{m_1(m_1-1)}{8}|\mathcal{H}|^2\psi\nonumber\\
&{}&+(l-2)\cdot\frac{m_1-1}{4}\mathcal{H}\cdot D_\beta\psi+\frac{m_1-1}{m_1}\left(D_1-\sum_{b=2}^l D_b\right)(D_\beta\psi)\nonumber\\
&{}&+(l-2)\cdot\frac{m_1-1}{2}\mathcal{H}\cdot\sum_{b=2}^l D_b\psi\nonumber\\
&=&\frac{m_1-1}{2}\sum_{a=1}^lD_a^\bot\mathcal{H}\cdot\psi+\frac{m_1(m_1-1)}{4}|\mathcal{H}|^2\psi\nonumber\\
&{}&+\frac{m_1-1}{m_1}\left(D_1-\sum_{b=2}^l D_b\right)(D_\beta\psi),
\eeq
where in the last step we used the relation
\beQ
\sum_{b=2}^lD_b\psi=\frac{1}{2}(D\psi-D_\beta\psi)
=\frac{m_1}{4}\mathcal{H}\cdot\psi-\frac{1}{2}D_\beta\psi.
\eeQ
Note, for $\alpha\in (I_1\cup I_2\cup\cdots\cup I_l)^c,$  we have
$\nabla_\alpha\psi=0.$ Thus
\beq
\left(D_1-\sum_{b=2}^lD_b\right)(D_\beta\psi)
=D_\beta^2\psi
&=&D^2\psi\nonumber\\
&=&D\left(\frac{m_1 }{2} \mathcal{H}\cdot\psi\right)\nonumber\\
&=&\frac{m_1}{2}(D^\bot\mathcal{H}\cdot\psi-\mathcal{H}\cdot D\psi)\nonumber\\
&=&\frac{m_1}{2}D^\bot\mathcal{H}\cdot\psi+\frac{m_1^2}{4}|\mathcal{H}|^2\psi.
\eeq
Therefore,
\beq
\frac{1}{2}R^M\psi
=\frac{m_1-1}{2}\sum_{a=1}^lD_a^\bot\mathcal{H}\cdot\psi+\frac{m_1-1}{2}D^\bot\mathcal{H}\cdot\psi
+\frac{m_1(m_1-1)}{2}|\mathcal{H}|^2\psi.
\eeq
This implies that $R^M=m_1(m_1-1)|\mathcal{H}|^2.$ The whole proof of Theorem 1 is complete.
\qed

Furthermore, now assume that
$m_1(R^M+\mathfrak{R}^N_{\psi})>(m_1-1)|\vec{H}|^2$ on $M_{\psi},$ the complement of which in $M^m$ is of zero measure, then one can choose $q$ such that
\beq\label{q}
(1-m_1q)^2=\frac{(m_1-1)|\vec{H}|}{\sqrt{\frac{m_1}{m_1-1}(R^M+\mathfrak{R}^N_{\psi})}-|\vec{H}|}\ \  \textup{on}\ \ M_{\psi},
\eeq
in (\ref{final}), to obtain
\beq
0&{\leq}&\int_M\left(1-\sqrt{\frac{m_2}{m_1}}\right)|T\psi|^2
+\left(\frac{1}{m_2}-\frac{1}{\sqrt{m_1m_2}}\right)|D_2\psi|^2\nonumber\\
&{}&+\cdots+\left(\frac{1}{m_k}
-\frac{1}{\sqrt{m_1 m_2}}\right)|D_k\psi|^2\nonumber\\
&\leq&\int_M(m_1q^2-2q+1)\left[\lambda_H^2
-\frac{1}{4}
\left(\sqrt{\frac{m_1}{m_1-1}(R^M+\mathfrak{R}^N_{\psi})}-|\vec{H}|\right)^2\right]|\psi|^2.
\quad\quad\quad
\eeq
At the same time by (\ref{q}), (\ref{Hcdot}) can be simplified as
\beq
\vec{H}\cdot\psi=\textup{sign}(\lambda_H)|\vec{H}|\omega_\bot\cdot\psi
\eeq
 in the limiting case.

This implies the following theorem
\begin{theorem}\label{theorem2}
Let $M^m\subset \overline{M}^{m+n}$ be a closed spin  submanifold of a Riemannian spin manifold $\overline{M}^{m+n}.$ Suppose $TM^m=T_1\oplus\cdots\oplus T_k,$ where $T_a$ are parallel and pairwise orthogonal distributions of dimension $m_a, a=1,\cdots,
k, $ and
$m_1\geq m_2\geq\cdots\geq m_k\geq0.$ Consider a non-trivial spinor field $\psi\in \Gamma(\mathbb{S})$ such that
$D_H\psi=\lambda_H\psi,$  and also assume that
$m_1(R^M+\mathfrak{R}^N_{\psi})>(m_1-1)|\vec{H}|^2$ on $M_{\psi},$  then
\beq\label{main result H2}
\lambda_H^2\geq\frac{1}{4}\min_{M_\psi}\left(\sqrt{\frac{m_1}{m_1-1}(R^M+\mathfrak{R}^N_{\psi})}-|\vec{H}|\right)^2.
\eeq
If $\lambda_H^2$ achieves its minimum and the normal bundle is flat,  the scalar curvature $R^M$ and $|\vec{H}|$ are both constants.
\end{theorem}

\begin{remark}  \ \\
$(1)$ If  $M^m\hookrightarrow \overline{M}^{m+1}$ is  not minimal, then the normal bundle  is  an oriented real line bundle, hence trivial. $(\ref{Hcdot})$ becomes into
\beq
(m_1-1)H=2(1-m_1q)^2\lambda_H.
\eeq
Therefore, we also get $\textup{sign}(\lambda_H)=\textup{sign}(H) $ in the limiting cases of  Theorem $1$ and Theorem $2.$\\
$(2)$ If $m(R^M+\mathfrak{R}^N_{\psi})>(m-1)|\vec{H}|^2$ on $M_{\psi},$ we just take $\beta=\textup{Id}$ in the proof of Theorem $2,$ then the result was obtained in \cite{Zhang, Hijazi1,Hijazi2, Morel, Ginoux}.
\\
$(3)$ If $M^m\hookrightarrow \overline{M}^{m+1}$ is a minimal closed spin hypersurface, the estimate $(\ref{main result H2})$  reduces to the Alexandrov's result \cite{A1}. Moreover by $(\ref{i}),(\ref{s})$ and $(\ref{alpha}),$ following the arguments in \cite{A1}, one knows that, if  $(\ref{main result H2})$ is  an
 equality,  the universal cover ${\tilde{M}}$ of $M$
is isometric to a product $M_1\times \cdots\times M_k,$ where $\textup{dim} M_s=m_s, s=1,\cdots,k, $
$ M_1$ has a real Killing
spinor and  $M_s$ has a parallel spinor if $m_s<m_1,$  and $M_s$ has a parallel spinor or a real
 Killing spinor if $m_s=m_1.$
\end{remark}

Note our approach can be also used to recover the estimate in  \cite {Moroianu}.  First,
it is not difficult to check that $D_\beta$ is a formally self-adjoint elliptic operator with respect to $L^2$-product on closed manifold,
only assumption that
$\textup{div}{\beta}:=(\nabla_{e_i}\beta)(e_i)=0$ is needed.
So if $\theta$ is  a  harmonic vector field of unit length,  we can define the $\beta$-twist
Dirac operator
\beq
D_\beta\phi:=e_i\cdot\nabla_{\beta(e_i)}\phi= \beta(e_i)\cdot\nabla_{e_i}\phi,
\eeq
where $\beta(X):=X-2\theta(X)\theta^{\sharp},\ \textup{for}\  X\in \Gamma(TM).$
The spectrum of $D_\beta$ is still discrete and real. Moreover, we also have
the following Lichnerowicz-type formula (see, \cite{Chen1})
\beq\label{Chenformula}
 D^2_\beta\phi
 &=&-\theta\cdot\nabla^*\nabla(\theta\cdot\phi)+\frac{R}{4}\phi.
\eeq

\quad

\begin{theorem}(\cite {Moroianu}, Moroianu and Ornea.)
Let $M^m\subset \tilde{M}^{m+1}$ be a minimal closed spin hypersurface of dimension $m\geq 3,$ and if on $M^m$ there exists a  harmonic 1-form $\theta$
of  unit length. Consider a non-trivial spinor field $\psi\in \Gamma(\mathbb{S})$ such that
$D\psi=\lambda\psi,$
then
\beq\label{estimate2}
\lambda^2\geq\frac{m-1}{4(m-2)}\min_MR^M.
\eeq
If $\lambda^2$ achieves its minimum,  $\theta$ is in fact a parallel 1-form and the eigenspinor $\psi$ corresponding to $\lambda$ satisfies that $D(\theta\cdot\psi)=\lambda\theta\cdot\psi.$
\end{theorem}
\proof In this case, we choose $q=\frac{1}{m-1}$
and
\begin{comment}
\beq
T_i\phi
&=&\nabla_i\phi
+\frac{\lambda}{2(m-1)}(\beta+\textup{Id})(e_i)\cdot e_0\cdot\phi+
\frac{1}{2}\nabla_{(\beta-\textrm{Id})(e_i)}\phi
\eeq
\end{comment}
for the first eigenspinor  $\psi$ for of $D,$ the min-max principle gives
\beq
\int_{M} |D_\beta\psi|^2
\stackrel{(\ref{Chenformula})}{=}\int_{M} |D(\theta\cdot\psi)|^2\geq \lambda^2\int_{M} |\psi|^2.
\eeq
Therefore, (\ref{key}) turns to be
\beQ
\mathfrak{R}e\int_M\left\langle \sum_{i=1}^me_i\cdot T_i\phi, \frac{1}{2}(D\phi-D_\beta\phi)\right\rangle
&\leq&\frac{\lambda}{2}\int_M\mathfrak{R}e\langle e_0\cdot \phi, D_\beta\phi-D\phi\rangle.
\eeQ
Hence, the same argument still works and the estimate (\ref{estimate2}) can be also obtained and in the limiting case,
$\nabla_\theta\psi=0,$
 and $T\psi=0,$ that is, for $i=1,\cdots,m-1,$
\beq
\nabla_i\psi=-\frac{\lambda}{m-1} e_i\cdot e_0\cdot \psi.
\eeq
Therefore, using the  $(\frac{1}{2}Ricci)$-formula yields
\beQ
\frac{1}{2}Ric(\theta)\cdot \psi
&=&D(\nabla_\theta\psi)-\nabla_\theta (D\psi)-\sum_1^{m}e_i\cdot \nabla_{\nabla_i\theta}\psi\\
&=&-\sum_1^{m}e_i\cdot \nabla_{\nabla_i\theta}\psi\\
&=&\frac{\lambda}{m-1}D\theta\cdot e_0\cdot\psi=0.
\eeQ
Hence, $Ric(\theta)=0.$ Furthermore, one obtains that $\theta$ is in fact parallel by
 Bochner-Weitzenb\"ock formula on $1$-forms
 \beq\label{formula on 1-forms}
\Delta=\nabla^*\nabla+Ric,
\eeq and then we can apply Theorem 3.1 in \cite{A} to know that  the universal covering space  of $M^m$ is a Riemannian product of the form $M_1\times \mathbb{R},$ where $M_1$ admits a
real Killing spinor.
At last,  for any 1-form $\theta$ and any spinor $\phi, $
we always have
\beq
D(\theta\cdot\phi)=-\theta\cdot D\phi-2\nabla_{\theta}\phi+(D\theta)\cdot \phi.
\eeq
\qed

\mysection{\textbf{Intrinsic estimates  for the twisted Dirac-Schr\"odinger operator}} \ls

In fact, the normal bundle of the submanifold can be replaced by an auxiliary arbitrary  vector bundle
on the submanifold and all the preceding computations can be done in an intrinsic way to obtain results for a twisted Dirac-Schr\"odinger operator on the manifold.

Let $(M^m,g )$ be  a  closed Riemannian spin manifold, and let $\Sigma M$ be the spinor bundle for $M^m$ with the canonical Riemannian connection $\nabla^{\Sigma M}.$ Let $E$ be any  vector bundle over $M^m $ equipped with
a metric  connection $\nabla^{E}$ and a spin structure.
Set
\beQ
\Sigma:=\Sigma M\otimes \Sigma E.
\eeQ
Recall that Clifford multiplication on $\Gamma(\Sigma)$ by a tangent vector field $X\in\Gamma(TM^m)$ is given by:
\beQ
X\cdot\phi:=(X\cdot\sigma)\otimes \varepsilon\in\Gamma(\Sigma),
\eeQ
for $\forall\ \phi=\sigma\otimes \varepsilon\in\Gamma(\Sigma)$ and the canonical tensor product connection on
$\Gamma(\Sigma)$ is defined by the formula
\beQ
\nabla(\sigma\otimes \varepsilon):=(\nabla^{\Sigma M}\sigma)\otimes \varepsilon+\sigma\otimes (\nabla^{E}\varepsilon),
\eeQ
where $\nabla^{\Sigma M}$ and $\nabla^{\Sigma E}$ denote the covariant derivatives  on $\Sigma M$ and
$\Sigma E$ respectively.
A direct verification shows that  the curvature transformation
$\mathcal{R}_{X,Y}=[\nabla_X,\nabla_Y]-\nabla_{[X,Y]}$
of $\Sigma M\otimes \Sigma E$
is also a derivation, i.e,
\beQ
\mathcal{R}(\sigma\otimes \varepsilon)
&=&(\mathcal{R}^M\sigma)\otimes \varepsilon+\sigma\otimes (\mathcal{R}^E\varepsilon),
\eeQ
where $\mathcal{R}^M$ and $ \mathcal{R}^E$ denote the curvature transformations of  $\Sigma M$ and
$\Sigma E$ respectively.

The Dirac operator on $M^m$ twisted with the bundle $\Sigma E$ is given by
\beQ
D^{\Sigma E}: \Gamma(\Sigma)&\longrightarrow &\Gamma(\Sigma)\\
\phi&\longmapsto & D^{\Sigma E}\phi=\sum_{i=1}^me_i\cdot\nabla_i\phi.
\eeQ
It is straightforward to check that  we have the following formula Schr\"odinger-Lichnerowicz-type formula
for the twisted Dirac operator $D^{\Sigma E}$ on the twisted spinor bundle $\Sigma=\Sigma M\otimes \Sigma E$
over $M^m$
\beq
(D^{\Sigma E})^2&=&\nabla^{*}\nabla+\frac{1}{2}\sum_{i,j=1}^me_i\cdot e_j\cdot\mathcal{R}_{e_i,e_j}\nonumber\\
&=&\nabla^{*}\nabla+\frac{1}{2}\sum_{i,j=1}^me_i\cdot e_j\cdot
(\mathcal{R}^M_{e_i,e_j}\otimes \textup{Id}+\textup{Id}\otimes \mathcal{R}^E_{e_i,e_j})\nonumber\\
&=&\nabla^{*}\nabla+\frac{1}{2}\sum_{i,j=1}^m\left[(e_i\cdot e_j\cdot
\mathcal{R}^M_{e_i,e_j})\otimes \textup{Id}+(e_i\cdot e_j\cdot\textup{Id})\otimes \mathcal{R}^E_{e_i,e_j}\right]\nonumber\\
&=&\nabla^{*}\nabla+\frac{1}{4}(R^M+\mathfrak{R}^E),
\eeq
where
$
\mathfrak{R}^E:\Sigma \rightarrow \Sigma
$
is a smooth symmetric bundle endomorphism defined by the formula
\beq
\mathfrak{R}^E(\sigma\otimes \varepsilon):=2\sum_{i,j=1}^m(e_i\cdot e_j\cdot \sigma)\otimes \mathcal{R}^E_{i,j}\varepsilon
\eeq
on vectors $\sigma\otimes \varepsilon\in \Gamma(\Sigma)$ of simple type.

For any smooth real function $F$ on $M^m,$ define the twisted Dirac-Schr\"odinger operator by
\beq
D_F=D^{\Sigma E}-\frac{1}{2}F.
\eeq
Now suppose $TM^m=T_1\oplus\cdots\oplus T_k,$ where $T_a$ are parallel and pairwise orthogonal
distributions of dimension $m_a, a=1,\cdots,
k, $ and
$m_1\geq m_2\geq\cdots\geq m_k.$
And  also assume
\beq
m_1(R^M+\kappa_1)>(m_1-1)F^2>0
\eeq
 on $M^m,$
 where $\kappa_1$ be the lowest eigenvalue of the endomorphism $\mathfrak{R}^E.$
Then define the modified connection on $\Sigma=\Sigma M\otimes \Sigma E$
\beq
T_i\phi
&=&\nabla_i\phi
+\left(\frac{1-q}{2(1-m_1q)}F+q\lambda_F\right) \frac{1}{2}(\beta(e_i)+e_i)\cdot\phi+
\frac{1}{2}\nabla_{(\beta-\textrm{Id})(e_i)}\phi,\quad
\eeq
where the smooth real function $q$ satisfies
\beq
(1-m_1q)^2=\frac{(m_1-1)|F|}{\sqrt{\frac{m_1}{m_1-1}(R^M+\kappa_1)}-|F|}.
\eeq
The same computations
as in the proof of Theorem \ref{theorem2}, lead to the following estimate,
\beq
&{}&\int_M\left(1-\sqrt{\frac{m_2}{m_1}}\right)|T\psi|^2
+\left(\frac{1}{m_2}-\frac{1}{\sqrt{m_1m_2}}\right)|D_2\psi|^2\nonumber\\
&{}&+\cdots+\left(\frac{1}{m_k}
-\frac{1}{\sqrt{m_1 m_2}}\right)|D_k\psi|^2\nonumber\\
&\leq&\int_M(m_1q^2-2q+1)\left[\lambda_F^2
-\frac{1}{4}\left(\sqrt{\frac{m_1}{m_1-1}(R^M+\kappa_1)}-|F|\right)^2\right]|\psi|^2,\quad
\eeq
where $D_F\psi=\lambda_F\psi$ for a non-trivial spinor field $\psi.$ Therefore we can conclude the following theorem.

\begin{theorem}
Let $(M^m,g)$ be a closed Riemannian spin manifold and $E\rightarrow M^m$ be a Riemannian and spin vector bundle of rank $n$ over $M^m.$ Assume
$TM^m=T_1\oplus\cdots\oplus T_k,$ where $T_a$ are parallel and pairwise orthogonal
distributions of dimension $m_a, a=1,\cdots,
k, $ and
$m_1\geq m_2\geq\cdots\geq m_k\geq0.$ Let  $\lambda_F$ be any
eigenvalue of the  Dirac-Schr\"odinger operator $D_F=D^{\Sigma E}-\frac{1}{2}F,$ associated with the eigenspinor $\psi.$ Assume $m_1(R^M+\kappa_1)>(m_1-1)F^2$ on $M^m,$ then
\beq\label{main result F}
\lambda_H^2\geq\frac{1}{4}\min_{M}\left(\sqrt{\frac{m_1}{m_1-1}(R^M+\kappa_1)}-|F|\right)^2.
\eeq
If $\lambda_F^2$ achieves its minimum,  $F$ are constant and also
$\mathfrak{R}^E\psi=\kappa_1\psi;$  In fact, the universal cover ${\tilde{M}}$ of $M^m$
is isometric to a product $M_1\times \cdots\times M_k,$ where $\textup{dim} M_s=m_s, M_1$ has a real Killing
spinor and  $M_s$ has a parallel spinor if $m_s<m_1,$  and $M_s$ has a parallel spinor or a real
 Killing spinor if $m_s=m_1.$
\end{theorem}
 Note there is a way to give an intrinsic meaning  to the modified connection (\ref{modicon})
 in the previous section, assuming an additional condition
 \beq
 \vec{H}\cdot\psi=F\omega_\bot\cdot\psi
 \eeq
 In fact, there exists an identification of   the restricted spinor bundle  $\mathbb{S}:={\Sigma \overline{M}|}_{M}$ with the twisted spinor bundle $\Sigma:=\Sigma M\otimes \Sigma NM$ or its direct sum $\Sigma\oplus\Sigma:$
 \beq
\mathbb{S}\longrightarrow
\left\{
    \begin{array}{cccccc}
\Sigma,\quad \textup{if}\ m\ \textup{or}\ n\ \textup{is even}\\
\Sigma\oplus\Sigma,\quad\quad\textup{otherwise},\\
\end{array}
    \right.
\eeq
which sends $\phi\in\Gamma(\mathbb{S})$ to $\phi^{*}\in \Gamma(\Sigma),$ for example, if $m$ or $n$ is even.
Moreover, with respect to this identification, Clifford multiplication by a vector field $X\in \Gamma(TM^m),$ is given (see \cite{Ginoux}) by
\beq
\forall \phi\in\Gamma(\mathbb{S}),\quad X\cdot\phi^{*}=(X\cdot\omega_{\bot}\cdot\phi)^{*}.
\eeq

\mysection{\textbf{Other generalized hypersurface Dirac operators}} \ls
In this section, we shall extend the above techniques to other generalized hypersurface Dirac operators which  appear in General Relativity. Let $p_{ij}$ be a-tensor on $\overline{M}^{m+1} $ and $P:=\bar{g}^{ij}p_{ij}|_M,$ then we
consider the following operator
\beq
D_\mathrm{P}=e_0\cdot D-\frac{\sqrt{-1}}{2}\mathrm{P}e_0\cdot
\eeq
Note that $D_\mathrm{P}$ is also formally self-adjoint with respect to $\int_M\langle\cdot,\cdot\rangle$ and, hence, $D_\mathrm{P}$ has real eigenvalues.
 Suppose $TM^m=T_1\oplus\cdots\oplus T_k,$ where $T_a$ are parallel and pairwise orthogonal distributions of dimension $m_a, a=1,\cdots,
k, $ and
$m_1\geq m_2\geq\cdots\geq m_k\geq0.$ We can  define
\beq
T_i\phi&=&\nabla_i+\frac{\sqrt{-1}}{4}p\mathrm{P}(\beta+\textup{Id})(e_i)\cdot\phi+\frac{1}{2}q\lambda_\mathrm{P}
e_0\cdot(\beta+\textup{Id})(e_i)\cdot\phi\nonumber\\
&{}&+\frac{1}{2}\nabla_{(\beta-\textrm{Id})(e_i)}\phi.
\eeq
Then for  $D_\mathrm{P}\psi=\lambda_\mathrm{P}\psi,$
\beq
|T\phi|^2
&=&|\nabla \psi|^2+\mathrm{P}(m_1pq-p-q)\mathfrak{R}e\langle D\psi,\sqrt{-1}\psi\rangle+(m_1q^2-2q)|D_\mathrm{P}\psi|^2\nonumber\\
&{}&+\frac{\mathrm{P}^2}{4}(m_1p^2+2q-2m_1pq)|\psi|^2-\sum_{s=m_1+1}^m|\nabla_s\psi|^2\nonumber\\
&{}&-\mathfrak{R}e\left\langle\left(-\frac{\sqrt{-1}}{2}p\mathrm{P}+q\lambda_\mathrm{P}e_0\cdot
\right)\psi,D\psi-D_\beta\psi\right\rangle
\eeq
In this case, we still have for any spinor field $\phi\in \Gamma(\mathbb{S})$
\beq\label{T2}
\sum_{i=1}^me_i\cdot T_i\phi=D\phi+m_1\left(-\frac{\sqrt{-1}}{2}p\mathrm{P}+q\lambda_\mathrm{P}e_0\cdot
\right) \phi+\frac{1}{2}(D_\beta\phi-D\phi).
\eeq
Hence
\beq
&{}&\mathfrak{R}e\int_M\langle \sum_{i=1}^me_i\cdot T_i\phi, -\frac{1}{2}(D_\beta\phi-D\phi)\rangle\nonumber\\
&{=}&\frac{m_1}{2}\mathfrak{R}e\int_M\left\langle \left(-\frac{\sqrt{-1}}{2}p\mathrm{P}+q\lambda_\mathrm{P}e_0\cdot
\right)\phi,D\phi-D_\beta\phi\right\rangle
\eeq
At the same time
\beQ
\int_M|D_\mathrm{P}\psi|^2
&=&\int_M|\nabla\phi|^2+\frac{1}{4} (R^M+\mathrm{P}^2)|\phi|^2-\mathrm{P}\mathfrak{R}e\langle D\phi,\sqrt{-1}\phi\rangle
\eeQ
Now assume
$m_1R^M>(m_1-1)\mathrm{P}^2$
on $M^m,$
 then
\beq\label{main result H3}
\lambda_\mathrm{P}^2\geq\frac{1}{4}\min_M\left(\sqrt{\frac{m_1}{m_1-1}R^M}-|\mathrm{P}|\right)^2.
\eeq


\begin{thebibliography} {99}



\bibitem{A1}B. Alexandrov, The first eigenvalue of the Dirac operator on locally reducible
           Riemannian manifolds, J. Geom. Phys. 57, 467--472(2007).

\bibitem{A} B. Alexandrov, G. Grantcharov, S. Ivanov, An estimate for the first eigenvalue
         of the Dirac operator on compact Riemannian spin manifold admitting a parallel
          one-form, J. Geom. Phys. 28, 263--270(1998).

\bibitem{Bar93}C. B\"ar, Real Killing spinors and holonomy, Comm. Math. Phys. 154, 509--521(1993).

\bibitem{Bar04}C. B\"ar, M. Dahl, The first Dirac eigenvalues on manifolds with positive
             scalar curvature, Proc. Amer. Math. Soc. 132, 3337--3344(2004).

\bibitem{Chen1}Y. Chen, The Dirac operator on manifold admitting parallel one-form, J. Geom. Phys.  117, 214--221 (2017).

\bibitem{Friedrich}T. Friedrich, Dirac operators in Riemannian geometry,
        Graduate Studies in Mathematics 25, American Mathematical Society(2000).

\bibitem{Gibbons} G. Gibbons, S. Hawking,  G. Horowitz,  M. Perry,  Positive mass theorems for black holes. Commun. Math. Phys. 88, 295--308 (1983).

\bibitem{Herzlich} M. Herzlich, The positive mass theorems for black holes revisited.
      J. Geom. Phys.  26, 97--111 (1998).


\bibitem{Ginoux} N. Ginoux, B. Morel, On eigenvalue estimate for the submanifold Dirac operator,
Int. J.  Math. 13, 533--548 (2002).


\bibitem{Hijazi1}O. Hijazi, X. Zhang, Lower bounds for the eigenvalues of Dirac operator, Part I. The hypersurface Dirac operators. Ann Glob Anal Geom. 19, 355--376(2001).

\bibitem{Hijazi2}O. Hijazi, X. Zhang, Lower bounds for the eigenvalues of Dirac operator, Part II. The submanifold Dirac operators. Ann Glob Anal Geom. 19, 163--181(2001).


\bibitem{K} E. C. Kim, Lower bounds of the Dirac eigenvalues on compact Riemannian spin manifolds
          with locally product structure, arXiv:math.DG/0402427 (2004).

\bibitem{Kirchberg86} K.-D. Kirchberg, An estimation for the first eigenvalue of the Dirac operator in closed K\"ahler manifolds of positive scalar curvature, Ann. Glob. Ann. Geom. 3, 291--325 (1986).

\bibitem{Kirchberg90} K.-D. Kirchberg, The first eigenvalue of the Dirac operator on K\"ahler manifolds, J. Geom. Phys. 7, 449--468(1990).

\bibitem{LM} H. B. Lawson,  M. L. Michelsohn, Spin Geometry,
         Princeton  Math Series. 38, Princeton University Press(1989).

\bibitem{Morel} B. Morel, Eigenvalue estimates for the  Dirac-Schr\"odinger operators,
               J. Geom. Phys. 38, 1--18(2001).

\bibitem{Moroianu} A. Moroianu, L. Ornea, Eigenvalue estimates for the Dirac operator and harmonic
              1-forms of constant length, C. R. Math. Acad. Sci. Paris. 338, 561--564 (2004).

\bibitem{Wang} Mc.K.Y. Wang, Parallel spinors and parallel forms, Ann. Glob. Anal. Geom. 7, 59--68(1989).

\bibitem{Yano} K. Yano, M. Kon,  Structures on manifolds. Singapore: World Sci.(1984).

\bibitem{Zhang} X. Zhang, Lower bounds for eigenvalues of hypersurface Dirac operators. Math Res Lett. 5, 199--210 (1998); A remark: Lower bounds for eigenvalues of hypersurface Dirac operators. Math Res Lett. 6, 465--466(1999).

\bibitem{Zhang} X. Zhang, Angular momentum and positive mass theorem, Commun. Math. Phys. 206, 137--155 (1999).
\end{thebibliography}
\end{document}